\documentclass[11pt]{amsart}
 
\usepackage[arrow,matrix,curve,cmtip,ps]{xy}
\usepackage{upref}

\newtheorem{thm}{Theorem}[section]   
\newtheorem{cor}[thm]{Corollary}     
\newtheorem{lem}[thm]{Lemma}         
\newtheorem{prop}[thm]{Proposition}  

\theoremstyle{definition}
\newtheorem{defn}[thm]{Definition}   

\theoremstyle{remark}
\newtheorem{rem}[thm]{Remark}        

\numberwithin{equation}{section}     

\newcommand{\secref}[1]{Section~\textup{\ref{#1}}}
\newcommand{\thmref}[1]{Theorem~\textup{\ref{#1}}}
\newcommand{\corref}[1]{Corollary~\textup{\ref{#1}}}
\newcommand{\lemref}[1]{Lemma~\textup{\ref{#1}}}
\newcommand{\propref}[1]{Proposition~\textup{\ref{#1}}}

\renewcommand{\c}[1]{\mathcal #1}

\renewcommand{\H}{\mathcal H}

\newcommand{\FF}{\mathcal F}

\newcommand{\A}{\mathcal A}
\newcommand{\D}{\mathcal D}

\newcommand{\UM}{UM}
\newcommand{\eps}{\epsilon}

\newcommand{\midtext}[1]{\quad\text{#1}\quad}
\newcommand{\righttext}[1]{\qquad\text{#1 }}

\DeclareMathOperator{\ad}{Ad}
\DeclareMathOperator{\id}{id}
\DeclareMathOperator{\im}{im}

\DeclareMathOperator{\ind}{Ind}
\DeclareMathOperator{\Ind}{Ind}
\DeclareMathOperator{\Rep}{Rep}

\DeclareMathOperator{\Aut}{Aut}

\DeclareMathOperator{\infl}{Inf}
\DeclareMathOperator{\Inf}{Inf}

\DeclareMathOperator*{\clsp}{\overline{span}}
\DeclareMathOperator*{\spn}{span}

\newcommand{\dec}{^{\text{dec}}}

\newcommand{\nonzero}[1]{\qquad\text{if $#1$
 \textup(and $0$ else\textup)}}

\newcommand{\lip}[1]{{}_{#1}\!}
\newcommand{\rip}[1]{_{\!#1}}

\newcommand{\<}{\langle}
\renewcommand{\>}{\rangle}

\newcommand{\pb}{q^*}
\newcommand{\what}{\widehat}
\newcommand{\lk}{\langle}
\newcommand{\rk}{\rangle}

\renewcommand{\:}{\colon}

\begin{document}

\title{Full duality for coactions of discrete groups}

\author{Siegfried Echterhoff}
\address{Westf\"alische Wilhelms-Universit\"at\\
SFB 478\\
Hittorfstr. 27\\
D-48149 M\"unster\\
Germany}
\email{echters@math.uni-muenster.de}

\author{John Quigg}
\address{Department of Mathematics\\
Arizona State University\\
Tempe, Arizona 85287}
\email{quigg@math.la.asu.edu}

\thanks{Research partially supported by National Science Foundation
Grant DMS9401253 and Deutsche Forschungsgemeinschaft (SFB 478)}


\date{June 19, 1999}

\begin{abstract}
Using the strong relation between coactions of a discrete group $G$
on $C^*$-algebras and Fell bundles over $G$ we prove a new version of
Mansfield's imprimitivity theorem for coactions of discrete groups.
Our imprimitivity theorem works for the universally defined full crossed
products and arbitrary subgroups of $G$ as opposed to the usual theory
of \cite{man, kq:imprimitivity} which uses the spatially defined reduced
crossed products and normal subgroups of $G$. Moreover, our theorem
factors through the usual one by passing to appropriate quotients. As
applications we show that a Fell bundle over a discrete group is amenable
in the sense of Exel \cite{exe:amenability} if and only if the double dual
action is amenable in the sense that the maximal and reduced crossed
products coincide. We also give a new characterization of induced
coactions in terms of their dual actions. 
\end{abstract}

\maketitle


\section{Introduction} 

One of the main tools in the study of crossed products by 
$C^*$-dynamical 
systems is Green's imprimitivity theorem (\cite[\S 2]{gre:local}): 
starting with an action
$\alpha:G\to \Aut(A)$ and a closed subgroup $H$ of $G$, it provides 
a Morita equivalence between the crossed products
$$C_0(G/H,A)\times_{\tau\otimes\alpha}G\quad\quad\text{and}
\quad\quad A\times_{\alpha|}H,$$
where $\tau$ denotes the translation action of $G$ on $G/H$.
Note that we can take either the full or the reduced
crossed products in the formulation of Green's theorem (see 
\cite[Theorem
3.15]{kas} and \cite[Lemma 4.1]{qs:regularity} for the reduced 
versions).
If $G$ is abelian, the imprimitivity theorem 
combined with Takesaki-Takai duality theory provides an even more
powerful tool for the investigation of the structure of certain 
crossed products. Note that
the connection of both theories is given by the existence of a 
natural 
isomorphism $$C_0(G/H,A)\times_{\tau\otimes\alpha}G\cong 
(A\times_{\alpha}G)\times_{\widehat{\alpha}|}\widehat{G/H},$$
where $\widehat{\alpha}$ denotes the dual action of the Pontrjagin
dual $\widehat{G}$ on $A\times_{\alpha}G$ (e.g., see 
\cite[\S7]{gre:local} and \cite{ech:indres, ech:trace} 
for details and applications).
If $G$ is nonabelian, Takesaki-Takai duality theory
becomes a duality between actions and coactions of $G$, so
in order to make these techniques available for the study 
of crossed products by actions and coactions of nonabelian groups it 
is necessary to have good working analogues of the imprimitivity 
theorem
for coactions of groups on 
$C^*$-algebras.

A major step towards a general imprimitivity theorem for coactions was
first achieved by Mansfield in \cite{man}, where he showed that,
if $\delta:A\to M(A\otimes C^*(G))$ is a coaction of $G$ on $A$ and
$N$ is a {\em normal} and {\em amenable} closed subgroup 
of $G$, there exists a natural Morita equivalence between
$$(A\times_{\delta}G)\times_{\widehat{\delta}|}N\quad\quad\text{and}
\quad\quad A\times_{\delta|}G/N,$$
where $\widehat{\delta}:G\to \Aut(A\times_{\delta}G)$ denotes the dual 
action of $G$ on $A\times_{\delta}G$.  This result was generalized to 
possibly non-amenable normal subgroups by the second author and 
Kaliszewski in \cite{kq:imprimitivity}, but with a technical 
restriction which is satisfied in particular when
 $\delta$ is a {\em 
normal} coaction in the sense of \cite[Definition 
2.1]{qui:fullreduced}.  For such coactions they obtained a natural 
Morita equivalence for the {\em reduced} crossed product 
$(A\times_{\delta}G)\times_{\widehat{\delta}|,r}N$ and $ 
A\times_{\delta|}G/N$.  Since every coaction $\delta$ has a {\em 
normalization} (see \cite[\S2]{qui:fullreduced}), this result should 
be interpreted as a general analogue of Green's theorem for reduced 
crossed products (note that, if $N$ is amenable, 
$B\times_{\beta}N\cong B\times_{\beta,r}N$ for every action $\beta$ of 
$N$, which explains why Mansfield's result can be regarded as a 
special case of the results in \cite{kq:imprimitivity}).

In \cite{ekr} it is shown that if we start with a dual coaction
$\widehat{\alpha}:A\times_{\alpha}G\to 
M\big((A\times_{\alpha}G)\otimes C^*(G)
\big)$ on the {\em full} crossed product $A\times_{\alpha}G$ 
by an action $\alpha:G\to\Aut(A)$, then there  exists also a full 
version of Mansfield's theorem.
More precisely, for $B=A\times_{\alpha}G$ and 
$\delta=\widehat{\alpha}$
the authors constructed a natural imprimitivity bimodule
between the full crossed products
$$B\times_{\delta}G\times_{\widehat\delta|}N
\quad\quad\text{and}\quad\quad B\times_{\delta|}G/N,$$
which actually factors through the reduced version 
of \cite{kq:imprimitivity} by passing to appropriate quotients.
Moreover, because
 $\delta$ is dual to an action,
the construction of the bimodule became much more natural and
it even allowed versions of Mansfield's imprimitivity theorems 
for possibly non-normal subgroups of $G$. 
Note that the big advantage of working with full crossed products 
rather than with the spatially defined reduced ones
is 
 that full crossed products can be defined via 
universal properties, 
which, in general,  are not available for the study of the reduced 
crossed products. 

The main purpose of this paper is to extend the results
of \cite{ekr}  to arbitrary coactions of {\em discrete} groups.
The main idea is to use the strong relationship  between
coactions of a discrete group
$G$ and Fell 
bundles over $G$ as observed by the second author in
\cite{qui:discrete} (see also \cite{eq:induced}): 
if $\delta \: A \to A \otimes  
C^*(G)$ is a coaction of the discrete group $G$ on the $C^*$-algebra  
$A$, then the set of spectral subspaces $\A = \{A_s : s \in G\}$ 
with  
$A_s = \{a \in A : \delta(a) = a \otimes s\}$ forms a Fell bundle 
over  
$G$.  For any such bundle we can form the {\em full} and the {\em  
reduced} cross-sectional algebras $C^*(\A)$ and $C_r^*(\A)$ of $\A$,  
which are completions of the $^*$-algebra $\Gamma_c(\A)$ of finitely  
supported sections with respect to the maximal and minimal  
topologically graded $C^*$-norms (see \cite{exe:amenability}).  We  
also have $\Gamma_c(\A) \subseteq A$ as a dense subalgebra, and the  
identity map on $\Gamma_c(\A)$ induces quotient maps 
$ C^*(\A) \to A \to C_r^*(\A)$.
There are canonical coactions $\delta^m \: C^*(\A) \to C^*(\A) 
\otimes  
C^*(G)$ and $\delta^n \: C_r^*(\A) \to C_r^*(\A) \otimes C^*(G)$  
determined by $a_s \mapsto a_s \otimes s$. The quotient maps  
$C^*(\A) \to A \to C_r^*(\A)$ are equivariant with respect to  
$\delta^m$, $\delta$, and $\delta^n$, and it is shown in \cite[Lemma  
2.1]{eq:induced} that
they induce isomorphisms 
\[ 
C^*(\A) \times_{\delta^m} G 
\xrightarrow{\cong} A \times_{\delta} G 
\xrightarrow{\cong} C^*_r(\A) \times_{\delta^n} G 
\]  
of the crossed products.  It is actually checked in  
\cite[Section 2]{eq:induced} that the coaction $\delta^n$ on the  
reduced cross-sectional algebra $C^*_r(\c A)$ constructed above  
coincides with the normalization of $\delta$. 
On the other side, $\delta^m$ should be regarded as a {\em 
maximalization}
of $\delta$, thus explaining our notation. 

As our main result, we
shall derive a new version of Mansfield's imprimitivity  
theorem, which works for full  
crossed products and arbitrary subgroups of a discrete group $G$.  
To be more precise,  
we show that for any subgroup $H$ of $G$ there exists a canonical  
$C^*(\A) \times_{\delta^m} G \times_{\widehat{\delta^m}} H - C^*(\A  
\times G/H)$ imprimitivity bimodule $X$, where $C^*(\A \times G/H)$  
denotes the maximal cross-sectional algebra of the Fell bundle $\A  
\times G/H$ over the transformation groupoid $G \times G/H$. We then 
show that this bimodule naturally  
factors to a $C_r^*(\A) \times_{\delta^n} G  
\times_{\widehat{\delta^n},r} H - C_r^*(\A \times G/H)$ 
imprimitivity  
bimodule for the reduced cross-sectional algebras on both sides.  
Note  
that if $H$ is normal, we have $C^*(\A \times G/H) = C^*(\A)  
\times_{\delta^m|} G/H$ and $C_r^*(\A \times G/H) = C_r^*(\A)  
\times_{\delta^n |} G/H$, where $\delta^m |$ and $\delta^n |$ denote  
the restrictions of $\delta^m$ and $\delta^n$ to $G/H$.  Using the  
fact that an arbitrary coaction $\delta \: A \to A \otimes C^*(G)$ 
of  
a discrete group $G$ ``lies between'' $\delta^m$ and $\delta^n$, we  
derive an ``intermediate'' imprimitivity theorem for $\delta$: if 
$H$  
is a normal subgroup of $G$ then $A \times_{\delta|} G/H$ is  
canonically Morita equivalent to a crossed product $A 
\times_{\delta}  
G \times_{\widehat{\delta},\mu} H$ ``lying between'' the full 
crossed  
product $A \times_{\delta} G \times_{\widehat{\delta}} H$ and the  
reduced crossed product $A \times_{\delta} G  
\times_{\widehat{\delta},r} H$.

We obtain an interesting consequence concerning  
amenability of Fell bundles in the sense of Exel  
\cite{exe:amenability}: a Fell bundle $\A$ over $G$ is amenable if 
and  
only if its dual action $\widehat{\delta^m}$ of $G$ on $C^*(\A)  
\times_{\delta^m} G$ is amenable in the (weak) sense that the full 
and  
reduced crossed products coincide (recall that Exel defined a Fell  
bundle $\A$ to be amenable if the regular representation $C^*(\A) 
\to  
C_r^*(\A)$ is an isomorphism).  This result potentially allows us to  
restrict questions related to amenability of Fell bundles to the  
special (and better understood) case of actions.  As a further  
application of our general imprimitivity theorem we shall derive a  
characterization of induced coactions as introduced in  
\cite{eq:induced}: a coaction $\delta^m \: C^*(\c A) \to C^*(\c A)  
\otimes C^*(G)$ is induced from a quotient group $G/H$ if and only 
if  
the dual action $\what{\delta^m}$ is twisted over $H$ in
 the sense 
of  
Green \cite{gre:local}. 


\section{Preliminaries} 
\label{prelim} 
 
Throughout \S\ref{prelim}--\S\ref{appl} $G$ will be a discrete  
group.  We adopt the conventions of \cite{eq:induced, 
qui:fullreduced,  
qui:discrete} for coactions of groups on $C^*$-algebras, and of
\cite{exe:amenability, fd} for Fell bundles.  We will need the more  
general notion of Fell bundles over discrete groupoids, for which we  
refer to \cite{kum:fell}.  If $\c B$ is a Fell bundle over a 
discrete  
groupoid $\c G$, we let $\Gamma_c(\c B)$ denote the $^*$-algebra of  
finitely supported sections.  A \emph{homomorphism} of $\c B$ into a  
$C^*$-algebra $C$ is a map $\phi \: \c B \to C$ which is linear on  
each fiber $B_s$, is as multiplicative as makes sense, and preserves  
adjoints.  A \emph{representation} of $\c B$ on a Hilbert space $\c 
H$  
is a homomorphism $\pi \: \c B \to \c L(\c H)$, and we say $\pi$ is  
\emph{nondegenerate} if $\clsp \{\pi(\c B) \c H\} = \c 
H$.  A  
$C^*$-algebra $B$ is the \emph{enveloping $C^*$-algebra} of a  
$^*$-algebra $B_0$ if the supremum of the $C^*$-seminorms on $B_0$ 
is  
finite and $B$ is the Hausdorff completion of $B_0$ in this largest  
$C^*$-seminorm. 
 
\begin{prop} 
\label{envelope} 
If $\c B$ is a Fell bundle over the discrete groupoid $\c G$, then 
the  
$^*$-algebra $\Gamma_c(\c B)$ has an enveloping $C^*$-algebra, which  
we denote by $C^*(\c B)$. 
\end{prop} 
 
\begin{proof} 
Let $\Pi$ be a representation of $\Gamma_c(\c B)$, and let $b =  
\sum_{x \in \c G} b_x \in \Gamma_c(\c B)$, with $b_x \in B_x$ 
finitely  
nonzero.  Using that $\Pi|_{B_{d(x)}}$ is a $*$-representation 
of the $C^*$-algebra $B_{d(x)}$ we
get 
$\|\Pi(b_x)\|^2=\|\Pi(b_x^*b_x)\|\leq
\|b_x^*b_x\|=\|b_x\|^2 $,
which implies that $\|\Pi(b)\|\leq\sum_x\|b_x\|$. But this 
suffices to establish the conclusion.
\end{proof}

\begin{cor} 
\label{rep} 
The assignment $\Pi \mapsto \Pi|_{\c B}$ gives a one-to-one  
correspondence between homomorphisms of $C^*(\c B)$ and of $\c B$, 
and  
between nondegenerate representations of $C^*(\c B)$ and of $\c B$. 
\end{cor} 
 
\begin{proof} 
Since there is an obvious one to one correspondence between 
homomorphisms
of $\c B$ and $\Gamma_c(\c B)$, the result follows from 
\propref{envelope}.
\end{proof} 
 
Let $\alpha$ be an action of $G$ on a $C^*$-algebra $B$, and let $B  
\times G$ be the associated semidirect product Fell bundle over $G$ 
as  
in \cite{fd}, so that $C^*(B \times G)$ is canonically isomorphic to  
the crossed product $B \times_\alpha G$.  For reference, the  
operations on the bundle are 
\[ (b,s)(c,t) = (b \alpha_s(c),st) \midtext{and} 
(b,s)^* = (\alpha_{s^{-1}}(b)^*,s^{-1}).  
\] 
Often it is convenient to restrict attention to a dense  
$^*$-subalgebra $B_0$ of $B$, and we write $\Gamma_c(B_0 \times G)$  
for the linear span of $B_0 \times G$ in $\Gamma_c(B \times G)$. 
 
\begin{lem} 
\label{action bounded} 
Let $(B, G, \alpha)$ be an action, and suppose $B$ is the enveloping  
$C^*$-algebra of an $\alpha$-invariant $^*$-subalgebra $B_0$.  Then 
$B  
\times_\alpha G$ is the enveloping $C^*$-algebra of $\Gamma_c(B_0  
\times G)$. 
\end{lem} 
 
\begin{proof} 
\propref{envelope} implies that $B\times_{\alpha}G$ is the enveloping 
$C^*$-algebra of $\Gamma_c(B\times G)$, so it suffices to show that
every representation $\Pi$ of $B_0\times G$ extends to a 
representation
of the Fell bundle $B\times G$. 
But using the representation $\pi$ of $B_0$ defined by $\pi(b)=\Pi(b,e)$ one 
can easily check that $\Pi$ restricts to a bounded linear map on each fiber 
$(B_0,s)$.
\end{proof} 
 
Now let $\delta$ be a coaction of $G$ on a $C^*$-algebra $A$, and 
let  
$\c A$ be the associated Fell bundle.  For any subgroup $H$ of $G$, 
we  
get a Fell bundle $\c A \times G/H$ over the discrete transformation  
groupoid $G \times G/H$.  For reference, the operations on $\c A 
\times G/H$ are 
\[ 
(a_s,trH)(b_t,rH) = (a_sb_t,rH) 
\midtext{and} 
(a_s,tH)^* = (a_s^*,stH). 
\] 
We define a notion of covariant representation for the ``restricted  
coaction'' $(A, G/H, \delta|)$ (whatever \emph{that} might be), and  
show in \propref{universal} that $C^*(\c A \times G/H)$ is  
characterized by a universal property for covariant 
representations.   
In \corref{cor-max} we deduce that $C^*(\c A \times G/H)$ reduces to  
the usual crossed product $A \times_{\delta|} G/H$ when $H$ is 
normal  
and $\delta$ is maximal (that is, $A = C^*(\c A)$).  Recall that
the {\em restriction} $\delta|$ of a coaction 
$\delta:A\to A\otimes C^*(G)$ to the quotient {\em group} $G/H$
is defined by 
$$\delta|=(\id_A\otimes q)\circ \delta:A\to A\otimes C^*(G/H),$$
where $q:C^*(G)\to C^*(G/H)$ denotes the quotient map.
\propref{universal} also shows that when $H$ is non-normal $C^*(\c A  
\times G/H)$ is analogous to
the ``full crossed product by the (dual) coaction of  
the homogeneous space $G/H$'' in \cite[Section 2]{ekr}. 
 
\begin{defn} 
A \emph{covariant representation of $(A, G/H, \delta|)$} is a pair  
$(\pi,\mu)$, where $\pi$ is a nondegenerate representation of $A$ 
and  
$\mu$ is a nondegenerate representation of $c_0(G/H)$ on the same  
Hilbert space, such that 
\begin{equation} 
\label{cov} 
\pi(a_s) \mu(\chi_{tH}) = \mu(\chi_{stH}) \pi(a_s) 
\righttext{for all} s,t \in G, a_s \in A_s, 
\end{equation} 
where $\chi_{tH}$ denotes the characteristic function of the 
singleton  
subset $\{tH\}$ of $G/H$. 
\end{defn} 
 
\begin{rem} 
When $H$ is normal, if follows from \cite[Lemma 2.2]{qui:discrete}  
that the above definition is equivalent to the more usual definition  
of covariant representation of the coaction $(A, G/H, \delta|)$,  
namely 
\[ 
(\pi \otimes \id) \circ \delta(a) 
= \ad \mu \otimes \id(w_{G/H})(\pi(a) \otimes 1) 
\righttext{for} a \in A, 
\] 
where $w_{G/H}$ is the unitary element of $M(c_0(G/H) \otimes  
C^*(G/H)) = c_b(G/H, C^*(G/H))$ defined by $w_{G/H}(sH) = sH$,
where on the right hand side we view $sH$ as an element of $C^*(G/H)$.  
However, when $H$ is non-normal there is no  
such object $w_{G/H}$. 
\end{rem} 
 
\begin{prop} 
\label{prop-rest} 
Let $(\pi,\mu)$ be a covariant representation of $(A, G, \delta)$, 
and  
let $\mu|$ denote the restriction $\mu|_{c_0(G/H)}$.  Then  
$(\pi,\mu|)$ is a covariant representation of $(A, G/H, \delta|)$. 
\end{prop} 
 
\begin{proof} 
We need only verify the covariance condition \eqref{cov}.  Take $s,t  
\in G$ and $a_s \in A_s$.  Since $\chi_{tH} = \sum_{h \in H}  
\chi_{th}$, the sum converging in the strict topology of $c_b(G) =  
M(c_0(G))$, we have 
\begin{align*} 
\pi(a_s) \mu|(\chi_{tH}) 
& = \pi(a_s) \sum_h \mu(\chi_{th}) 
= \sum_h \pi(a_s) \mu(\chi_{th}) 
= \sum_h \mu(\chi_{sth}) \pi(a_s) 
\\& = \biggl(\sum_h \mu(\chi_{sth})\biggr) \pi(a_s) 
= \mu|(\chi_{stH}) \pi(a_s). 
\end{align*} 
\end{proof} 
 
\begin{prop} 
\label{universal} 
For every covariant representation $(\pi,\mu)$ of $(A, G/H, 
\delta|)$,  
there exists a unique nondegenerate representation $\pi \times \mu$ 
of  
$C^*(\c A \times G/H)$ such that 
\begin{equation} 
\label{int} 
\pi \times \mu(a,tH) = \pi(a) \mu(\chi_{tH}) 
\righttext{for all} a \in \c A, t \in G. 
\end{equation} 
Moreover, in the case $A = C^*(\c A)$, every nondegenerate  
representation of $C^*(\c A \times G/H)$ arises this way. 
\end{prop} 
 
\begin{proof} 
First assume $(\pi,\mu)$ is a covariant representation of $(A, G/H,  
\delta|)$.  It is easy to check that $(a,tH) \mapsto \pi(a)  
\mu(\chi_{tH})$ is a representation of the Fell bundle $\c A \times  
G/H$, and moreover this representation is nondegenerate since both  
$\pi$ and $\mu$ are.  So, by \corref{rep} there exists a unique  
representation $\pi \times \mu$ of $C^*(\c A \times G/H)$ satisfying  
the compatibility condition \eqref{int}. 
 
Conversely, assume $A = C^*(\c A)$, and let $\Pi$ be a nondegenerate  
representation of $C^*(\c A \times G/H)$.  For each coset $tH \in 
G/H$  
define a representation $\sigma_{tH}$ of the unit fiber algebra 
$A_e$  
by 
\[ 
\sigma_{tH}(a_e) = \Pi(a_e,tH). 
\] 
Fix a bounded approximate identity $\{d_i\}$ for $A_e$, and for each  
$tH \in G/H$ put 
\[ 
p_{tH} = \lim \sigma_{tH}(d_i), 
\] 
the limit taken in the weak operator topology.  Then $\{p_{tH} : tH  
\in G/H\}$ is an orthogonal family of projections, hence determines 
a  
representation $\mu$ of $c_0(G/H)$ such that $\mu(\chi_{tH}) =  
p_{tH}$.  Computations similar to those in the proof of 
\cite[Theorem  
3.3]{eq:induced} show that for all $a \in \c A$ the sum 
\[ 
\pi(a): = \sum_{tH \in G/H} \Pi(a,tH) 
\] 
converges in the weak operator topology and determines a  
representation $\pi$ of $A$, and moreover $(\pi,\mu)$ is a covariant  
representation of $(A, G/H, \delta|)$ such that $\Pi = \pi \times  
\mu$. 
\end{proof} 
 
Notice that if $A\neq C^*(\A)$, the case 
$H=G$ shows that there might be strictly
fewer covariant representations of $(A,G/H,\delta|)$ than
representations of $C^*(\A\times G/H)$. From 
\propref{universal} and the universal properties of
crossed products by coactions with respect to covariant
representations we get

\begin{cor} 
\label{cor-max} 
If $H$ is normal then there is a unique isomorphism $C^*(\c A) 
\times_{\delta^m|} G/H \cong C^*(\c A \times G/H)$ mapping $j_{C^*(\c 
A)}(a_s) j_{G/H}(\chi_{tH})$ to $(a_s,tH)$
\textup(where $j_{C^*(\c A)}$ and $j_{G/H}$ denote the canonical maps 
of $C^*(\c A)$ and $c_0(G/H)$ into 
$M(C^*(\A)\times_{\delta^m|}G/H)$\textup).
\end{cor} 
 
Now we recall from \cite{kum:fell} Kumjian's construction of what we  
call the \emph{reduced} cross-sectional algebra of a Fell bundle $\c  
B$ over a groupoid $\c G$.  Kumjian only needed $\c G$ to be  
$r$-discrete, but our $\c G$ will actually be discrete.  
$\Gamma_c(\c  
B)$ is given a pre-Hilbert $\Gamma_c(\c B^0)$-module structure 
(where  
$\c B^0$ is the restricted bundle $\c B|_{\c G^0}$ over the unit 
space  
$\c G^0$), with inner product 
\[ 
\< f,g \>\rip{\Gamma_c(\c B^0)} = (f^*g)|_{\c G^0} 
\righttext{for} f,g \in \Gamma_c(\c B). 
\] 
Here $\Gamma_c(\c B^0)$ is regarded as a subalgebra of $\Gamma_0(\c  
B^0)$, the $C_0$-section algebra of the $C^*$-bundle $\c B^0$.  Then  
the completion $L^2(\c B)$ of $\Gamma_c(\c B)$ is a full Hilbert  
$\Gamma_0(\c B^0)$-module.  Left multiplication in $\Gamma_c(\c B)$  
extends to a nondegenerate action of $\Gamma_c(\c B)$ on $L^2(\c B)$  
by adjointable operators, and the completion of $\Gamma_c(\c B)$ in  
the norm of this representation is the reduced cross-sectional 
algebra  
$C^*_r(\c B)$.  The identity map on $\c B$ extends uniquely to a  
surjective homomorphism, which we call the \emph{regular  
representation}, of $C^*(\c B)$ onto $C^*_r(\c B)$.  Further, we 
call  
a Fell bundle $\c B$ over a discrete groupoid $\c G$ \emph{amenable}  
if the regular representation of $C^*(\c B)$ onto $C^*_r(\c B)$ is  
faithful.  Note that in the special case where $\c G$ is actually a  
group, this is the same terminology used by Exel  
\cite{exe:amenability}.  Note also that in the case of an action 
$(B,  
G, \alpha)$, $C^*_r(B \times G)$ can be identified with the reduced  
crossed product $B \times_{\alpha,r} G$. 
 
Consider the case $\c B = \c A \times G/H$, where $\c A$ is a Fell  
bundle over $G$, and $\c G = G \times G/H$.  The unit space $\c G^0 
=  
\{e\} \times G/H$ can be identified with $G/H$.  Then $\c B^0 = A_e  
\times G/H$, a trivial $C^*$-bundle, so the $C_0$-section algebra is 
\[ 
\Gamma_0(A_e \times G/H) = A_e \otimes c_0(G/H) = c_0(G/H,A_e). 
\] 
The inner product on $\Gamma_c(\c B) = \Gamma_c(\c A \times G/H)$ is 
\[ 
\< f,g \>\rip{A_e \otimes c_0(G/H)} = (f^*g)|_{\{e\} \times G/H}. 
\] 
Writing $f$ as a finitely nonzero sum $\sum_{s,tH}(f_{s,tH},tH)$, 
and  
similarly for $g$, we arrive at 
\begin{align*} 
\biggl\lk \sum_{s,tH} (f_{s,tH},tH), 
\sum_{u,vH} (g_{u,vH},vH) \biggr\rk 
\rip{A_e \otimes c_0(G/H)} 
& = \sum_{s,tH} (f^*_{s,tH} g_{s,tH},tH) 
\\& = \sum_{tH \in G/H} 
\biggl(\sum_{s \in G} \< f_{s,tH},g_{s,tH} \>\rip{A_e},tH\biggr). 
\end{align*}
In particular, for generators $(a_s, tH), (b_u, vH) \in \c A \times  
G/H$ we have 
\[ 
\< (a_s, tH), (b_u, vH) \>\rip{A_e \otimes c_0(G/H)} 
= (a_s^* b_u, tH) 
\nonzero{s = u, tH = vH}. 
\] 
The regular representation, which we denote by
 $\Lambda$,
of the cross-sectional algebra $C^*(\c A \times G/H)$ on the Hilbert  
$A_e \otimes c_0(G/H)$-module 
$L^2(\c A \times G/H)$ is given by
\[ 
\Lambda(a_s, tH) (b_u \otimes \chi_{vH}) 
= a_s b_u \otimes \chi_{vH} 
\nonzero{tH = vH}, 
\] 
and $C^*_r(\c A \times G/H)$ is the image of $C^*(\c A \times G/H)$ 
in  
$\c L(L^2(\c A \times G/H))$. 
 
We want to identify the reduced cross-sectional algebra $C^*_r(\c A  
\times G/H)$ with the image $\im j_{C^*(\c A)} \times j_G|$ of 
$C^*(\A\times G/H)$ in $M(C^*(\c A) \times_{\delta^m} G)$,
where $j_{C^*(\c A)}$ and $j_G$ denote the canonical maps of 
$C^*(\c A)$ and $c_0(G)$ into $M(C^*(\c A) \times_{\delta^m} G)$
and $j_{C^*(\c A)} \times j_G|$ is the corresponding $*$-homomorphism
of $C^*(\c A\times G/H)$ as described in Proposition \ref{universal}.
For this we need the following technical lemma on Hilbert modules, 
which is an  
easy modification of \cite[Lemma 2.5]{sie:morita}.  
The essential idea is that a  
linear map between Hilbert modules which preserves inner products is  
automatically a module homomorphism.  Let $C$ and $D$ be  
$C^*$-algebras, let $Z$ be a (right) Hilbert $D$-module, and suppose  
$C$ is represented by adjointable operators on $Z$.  If $Z$ is full 
as  
a Hilbert $D$-module and the action of $C$ on $Z$ is nondegenerate, 
we  
say $Z$ is a \emph{right-Hilbert $C - D$ bimodule}.  We use the  
notation ${}_CZ_D$ to indicate that the coefficient algebras of the  
bimodule $Z$ are $C$ and $D$.  If $C_0$ and $D_0$ are dense  
$^*$-subalgebras and $Z_0$ is a dense linear subspace such that $C_0  
Z_0 \cup Z_0 D_0 \subseteq Z_0$ and $\< Z_0, Z_0 \>\rip D \subseteq  
D_0$, we say ${}_CZ_D$ is the \emph{completion} of the  
\emph{right-pre-Hilbert bimodule} ${}_{C_0}(Z_0)_{D_0}$. 
 
\begin{lem}[{cf.\ \cite[Lemma 2.5]{sie:morita}}]
\label{module hom}
Suppose ${}_CZ_D$ and ${}_EW_F$ are right-Hilbert bimodules such 
that  
${}_CZ_D$ is the completion of a right-pre-Hilbert bimodule  
${}_{C_0}(Z_0)_{D_0}$, and suppose we are given homomorphisms $\phi 
\:  
C \to E$ and $\psi \: D \to F$ and a linear map $\Psi \: Z_0\to W$  
with dense range such that for all $c \in C_0$ and $z,w \in Z_0$ we  
have 
\begin{enumerate} 
\item 
$\Psi(cz) = \phi(c) \Psi(z)$, and 
 
\item 
$\< \Psi(z),\Psi(w) \>\rip{F} = \psi\bigl(\< z,w \>\rip{D}\bigr)$. 
\end{enumerate} 
Then $\Psi$ extends uniquely to a right-Hilbert bimodule 
homomorphism  
of ${}_CZ_D$ onto ${}_EW_F$.  Moreover, if $Z$ and $W$ are actually 
imprimitivity bimodules, this extension of $\Psi$ is an 
imprimitivity  
bimodule homomorphism. 
\end{lem}

\begin{prop}\label{prop-reduced} 
With the notation preceding \lemref{module hom}, the two 
homomorphisms  
$\Lambda$ and $j_{C^*(\c A)} \times j_G|$ of $C^*(\c A \times G/H)$  
have the same kernel. 
\end{prop} 
 
\begin{proof} 
It suffices to produce isomorphic right-Hilbert $C^*(\c A \times G/H) 
- A_e$ bimodules $Y$ and $Z$ such that the homomorphisms of $C^*(\c  
A \times G/H)$ into $\c L(Y)$ and $\c L(Z)$ have the same kernels as  
$\Lambda$ and $j_{C^*(\c A)} \times j_G|$, respectively.  To get 
$Y$,  
note that $A_e \otimes c_0(G/H)$ acts faithfully by adjointable  
operators on the external tensor product $A_e \otimes \ell^2(G)$ of  
the Hilbert $A_e$-module $A_e$ and the Hilbert space $\ell^2(G)$, via 
\[ 
(a \otimes f) (b \otimes g) = ab \otimes fg 
\righttext{for} a,b \in A_e, f \in c_0(G/H), g \in \ell^2(G). 
\] 
Thus, the associated homomorphism of $C^*(\c A \times G/H)$ into the  
adjointable operators on the balanced tensor product $L^2(\c A 
\times  
G/H) \otimes_{A_e \otimes c_0(G/H)} (A_e \otimes \ell^2(G))$ has the  
same kernel as $\Lambda$.  

In order to construct the bimodule $Z$, first note that,
since the coaction $\delta^n$ on $C^*_r(\c A)$ is the normalization 
of $\delta^m$,
the
crossed product $C^*(\A)\times_{\delta^m}G$ has a faithful 
representation
on the Hilbert $A_e$-module
 $L^2(\c A) \c \otimes \ell^2(G)$.
The action is given on the generators by
\[
(a_s, t) (b_u \otimes \chi_v)
= a_s b_u \otimes \lambda_s M_{\chi_t} \chi_v
= a_s b_u \otimes \chi_{sv}
\nonzero{t = v}.
\]
If we compose this representation with 
$j_{C^*(\c A)} \times j_G|:C^*(\A\times G/H)\to 
M(C^*(\A)\times_{\delta^m}G)$,
$C^*(\c A \times G/H)$ acts on $L^2(\c A) \times \ell^2(G)$ by
\[
(a_s, tH) (b_u \otimes \chi_v)
= a_s b_u \otimes \chi_{sv}
\nonzero{v \in tH}.
\]
Thus for $Z$
we take $L^2(\c A) \otimes  \ell^2(G)$ equipped with this action. 
 
We want to define an isomorphism $\Psi$ between the Hilbert  
$A_e$-modules $L^2(\c A \times G/H) \otimes_{A_e \otimes c_0(G/H)}  
(A_e \otimes \ell^2(G))$ and $L^2(\c A) \otimes \ell^2(G)$.  We 
begin  
by defining $\Psi$ on the generators: 
\[ 
\Psi((a_s, tH) \otimes (b \otimes \chi_r)) 
= a_s b \otimes \chi_{sr} 
\nonzero{r \in tH}. 
\] 
and then extending additively.  This gives a linear map from  
$\Gamma_c(\c A \times G/H) \odot (A_e \odot C_c(G))$ to $L^2(\c A)  
\otimes \ell^2(G)$ with dense range.  We show $\Psi$ preserves inner  
products.  For $a_s, b_u \in \c A$, $tH, vH \in G/H$, $c, d \in 
A_e$,  
and $r, w \in G$ we have 
\begin{align*} 
& \< \Psi\bigl( (a_s, tH) \otimes (c \otimes \chi_r) \bigr), 
\Psi\bigl( (b_u, vH) \otimes (d \otimes \chi_w) \bigr) \>\rip{A_e} 
\\ & \quad =  
\< a_s c \otimes \chi_{sr}, b_u d \otimes \chi_{uw} \>\rip{A_e} 
\nonzero{r \in tH, w \in vH} 
\\ & \quad =  
\< a_s c, b_u d \>\rip{A_e} 
\< \chi_{sr}, \chi_{uw} \> 
\nonzero{r \in tH, w \in vH} 
\\ & \quad =  
c^* a^*_s b_u d 
\nonzero{ s = u, sr = uw , r \in tH, w \in vH}. 
\end{align*} 
On the other hand, 
\begin{align*} 
& \< (a_s, tH) \otimes (c \otimes \chi_r), 
(b_u, vH) \otimes (d \otimes \chi_w) \>\rip{A_e} 
\\ & \quad =  
\< c \otimes \chi_r,  
\< (a_s, tH), (b_u, vH) \>\rip{A_e \otimes c_0(G/H)} 
(d \otimes \chi_w) \>\rip{A_e} 
\\&\quad =  
\< c \otimes \chi_r, (a^*_s b_u, tH) (d \otimes \chi_w) \>\rip{A_e} 
\nonzero{s = u, tH = vH} 
\\&\quad =  
\< c \otimes \chi_r, a^*_s b_u d \otimes \chi_w \>\rip{A_e} 
\nonzero{s = u, tH = vH, w \in tH} 
\\&\quad =  
c^* a^*_s b_u d \< \chi_r, \chi_w \> 
\nonzero{s = u, tH = vH, w \in tH} 
\\&\quad =  
c^* a^*_s b_u d 
\nonzero{s = u, tH = vH, w \in tH, r = w}. 
\end{align*} 
Since the two sets of conditions 
\[ 
\{ s = u, sr = uw , r \in tH, w \in vH \} 
\midtext{and} 
\{ s = u, tH = vH, w \in tH, r = w\} 
\] 
are equivalent, $\Psi$ preserves the inner products. 
 
A similar computation shows that $\Psi$ preserves the 
left module actions. Thus, 
by \lemref{module hom}, $\Psi$ extends to a surjective right-Hilbert  
bimodule homomorphism, which we still denote by $\Psi$.  Since the  
right hand coefficient homomorphism is the identity map on $A_e$,  
$\Psi$ is actually an isomorphism.
\end{proof} 
 
\begin{rem} 
When $H = \{e\}$ the result shows that the reduced  
cross-sectional algebra $C^*_r(\c A \times G)$ of the Fell bundle 
$\c  
A \times G$ over $G \times G$ is isomorphic to the crossed product  
$C^*(\c A) \times_{\delta^m} G$.  Hence, the Fell bundle $\c A 
\times  
G$ is always amenable in the sense that the regular representation  
$\Lambda \: C^*(\c A \times G) \to C^*_r(\c A \times G)$ is faithful. 
\end{rem} 
 
\begin{cor} 
\label{cor-min} 
If $H$ is normal then there is a unique isomorphism $C^*_r(\c A)  
\times_{\delta^n|} G/H \cong C^*_r(\c A \times G/H)$ mapping  
$j_{C^*_r(\c A)}(a_s) j_{G/H}(\chi_{tH})$ to $(a_s,tH)$. 
\end{cor} 
 
\begin{proof} 
We must show that $j_{C^*_r(\c A)} \times j_{G/H}$ has the same kernel as  
the regular representation $\Lambda \: C^*(\c A \times G/H) \to  
C^*_r(\c A \times G/H)$. By \propref{prop-reduced}, it suffices to  
show that
the homomorphisms $j_{C^*_r(\c A)} \times j_{G/H}$ and  
$j_{C^*(\c A)} \times j_G|$ of $C^*(\c A \times G/H)$ have the same  
kernel. These maps fit into a commutative diagram 
\[ 
\xymatrix{ 
{C^*(\c A \times G/H)} 
\ar[rr]^{j_{C^*(\c A)} \times j_G|} 
\ar[d]_{j_{C^*_r(\c A)} \times j_{G/H}} 
&& {M(C^*(\c A) \times_{\delta^m} G)} 
\ar[d]^{j_{C^*_r(\c A)} \times j_G} 
\\ 
{C^*_r(\c A) \times_{\delta^n|} G/H} 
\ar[rr]_{j_{C^*_r(\c A)} \times j_G|} 
&& {M(C^*_r(\c A) \times_{\delta^n} G)} 
}. 
\] 
The right hand map is an isomorphism. Since the coaction $\delta^n$ 
is normal, by \cite[Lemma 3.1]{kq:imprimitivity} the bottom map is 
injective. The result follows.
\end{proof}


\section{The imprimitivity theorem} 
\label{main} 
 
In this section we prove the imprimitivity theorems for coactions of  
discrete groups.  Starting with a Fell bundle $\A$ over $G$, we 
first  
construct a $C^*(\c A) \times_{\delta^m} G \times_{\what{\delta^m}} 
H  
- C^*(\c A \times G/H)$ imprimitivity bimodule $X$, where 
${\delta^m}$  
is the canonical coaction of $G$ on $C^*(\A)$.  As usual, we work 
with  
dense subspaces.  For $C^*(\A \times G/H)$ we take the dense  
$^*$-subalgebra $C_0: = \Gamma_c(\A \times G/H)$, where we remind 
the  
reader to regard $\A \times G/H$ as a Fell bundle over the groupoid 
$G  
\times G/H$.  For $C^*(\A) \times_{\delta^m} G  
\times_{\widehat{\delta^m}} H$ we form the corresponding dense  
$^*$-subalgebra $B_0: = \Gamma_c(\A \times G \times H)$, and we put  
$X_0 = \Gamma_c(\A \times G)$.  For reference, the operations on the  
$^*$-algebra $B_0$ are given on the generators by 
\begin{align} 
\label{eq-b0} 
\begin{split} 
(a_s,t,h)(a_u,v,l)& = (a_sa_u,vh^{-1},hl) 
\nonzero{th = uv} \\ 
(a_s,t,h)^*& = (a_s^*,sth,h^{-1}). 
\end{split} 
\end{align} 
The $B_0 - C_0$ pre-imprimitivity bimodule structure on $X_0$ will 
be  
given on the generators by 
\begin{align}\label{eq-X0} 
\begin{split} 
(a_s,t) \cdot(a_u,vH)& = (a_sa_u,u^{-1}t) 
\nonzero{tH = uvH} \\ 
(a_q,r,h) \cdot(a_s,t)& = (a_qa_s,th^{-1}) 
\nonzero{rh = st} \\ 
\< (a_s,t),(a_u,v) \>\rip{C_0}& = (a_s^*a_u,vH) 
\nonzero{st = uv} \\ 
\lip{B_0}\< (a_s,t),(a_u,v) \>& = (a_sa_u^*,ut,t^{-1}v) 
\nonzero{tH = vH}. 
\end{split} 
\end{align} 
 
\begin{thm} 
\label{imp} 
Suppose that $\A$ is a Fell bundle over the discrete group $G$ and 
$H$  
is a subgroup of $G$.  The above operations make $X_0$ into a $B_0 -  
C_0$ pre-imprimitivity bimodule.  Consequently, its completion is a  
$C^*(\A) \times_{\delta^m} G \times_{\what{\delta^m}} H - C^*(\A  
\times G/H)$ imprimitivity bimodule $X$. 
\end{thm} 
 
\begin{proof} 
We closely follow the lines of the proof of \cite[Theorem  
4.1]{eq:induced}, where we proved an imprimitivity theorem for 
crossed  
products by induced coactions.  We have to check the following items: 
\begin{enumerate} 
\item 
$X_0$ is a $B_0 - C_0$ bimodule; 
 
\item $\lip{B_0}\< b \cdot x,y \> = b \lip{B_0}\< x,y \>$ 
and $\< x,y \cdot c \>\rip{C_0} = \< x,y \>\rip{C_0} c$; 
 
\item $\lip{B_0}\< x,y \>^* = \lip{B_0}\< y,x \>$ 
and $\< x,y \>\rip{C_0}^* = \< y,x \>\rip{C_0}$; 
 
\item $\lip{B_0}\< x,y \>$ is linear in $x$ 
and $\< x,y \>\rip{C_0}$ is linear in $y$; 
 
\item $x \cdot \< y,z \>\rip{C_0} = \lip{B_0}\< x,y \> \cdot z$; 
 
\item $\spn \lip{B_0}\< X_0,X_0 \>$ is dense in $B_0$ 
and $\spn \< X_0,X_0 \>\rip{C_0}$ is dense in $C_0$; 
 
\item $\lip{B_0}\< x,x \> \ge 0$ and $\< x,x \>\rip{C_0} \ge 0$; 
 
\item $\< b \cdot x, b \cdot x \>\rip{C_0} \le \|b\|^2 \< x,x  
\>\rip{C_0}$ and $\lip{B_0}\< x \cdot c, x \cdot c \> \le \|c\|^2  
\lip{B_0}\< x,x \>$. 
\end{enumerate} 
The verifications of the algebraic properties (i)--(v) are routine 
and  
we omit the details on this.  As in \cite{eq:induced}, we prove  
(vi)--(vii) in one whack using Rieffel's trick: it suffices to 
produce  
nets in both $B_0$ and $C_0$, each term of which is a finite sum of  
the form $\sum \< x_i,x_i \>$, which are approximate identities for  
both the algebras and the module multiplications in the inductive  
limit topologies. 
 
To construct an approximate identity for $C_0$ let $\{a_i\}_{i \in 
I}$  
be a bounded and positive approximate identity for the unit fiber  
$A_e$, and let $\FF$ denote the family of finite subsets of $G/H$,  
directed by inclusion.  For each $F \in \FF$ we choose $S_F 
\subseteq  
G$ comprising exactly one element from each coset in $F$.  Then 
\[ 
\sum_{t \in S_F} \< (a_i^{1/2},t), (a_i^{1/2},t) \>\rip{C_0} 
= \sum_{tH \in F} (a_i, tH) 
\] 
for all $(i,F) \in I \times \FF $, which implies 
\[ 
\left\{ \sum_{t \in S_F} \< (a_i^{1/2},t), (a_i^{1/2},t) \>\rip{C_0} 
\right\}_{(i,F) \in I \times \FF}  
\] 
is an approximate identity (in the inductive limit topologies) for  
both the algebra $C_0$ and the right module multiplication of $C_0$ 
on  
$X_0$.  For $B_0$, we let $\tilde{\FF}$ denote the finite subsets of  
$G$.  Then the desired approximate identity in $B_0$ is given by the  
elements 
\[ 
\sum_{t \in F} \lip{B_0}\< (a_i^{1/2},t), (a_i^{1/2},t) \> 
= \sum_{t \in F}(a_i,t,e), 
\] 
where $(i,F)$ runs through the directed set $I \times \tilde{\FF}$. 
 
In order to show (viii) we first observe that each generator  
$(a_s,t,h) \in B_0$ determines an adjointable operator on the  
pre-Hilbert $C_0$-module $X_0$ with adjoint given by the action of  
$(a_s,t,h)^* = (a_s^*, sth,h^{-1})$.  The product $(a_s,t,h)^*(a_s,  
t,h) = (a_s^*a_s, t,e)$ is a positive element of the $C^*$-algebra  
$(A_e, t, e)$.  Extend the action of $(A_e,t,e)$ on $X_0$ to its  
unitization $(\tilde{A}_e, t,e)$ in the obvious way.  We find 
elements  
$b_e \in \tilde{A}_e$ such that $(b_e^*b_e,t,e) = \|(a_s, t,  
h)\|^21-(a_s^*a_s,t,e)$.  Hence, using (vii), we get 
\[ 
\|(a_s,t,h)\|^2 \< x,x \>\rip{C_0} 
- \< (a_s,t,h) \cdot x,(a_s,t,h) \cdot x \>\rip{C_0} 
= \< (b_e,t,e) \cdot x,(b_e,t,e) \cdot x \>\rip{C_0}  
\geq 0. 
\] 
Thus the action of $(a_s,t,h)$ extends to an adjointable operator on  
the Hilbert $C^*(\c A \times G/H)$-module completion $X$, which by 
(i)  
gives us a $^*$-homomorphism of $B_0$ into the $C^*$-algebra of  
adjointable module maps on $X$.  By \lemref{action bounded} this  
homomorphism extends to the enveloping $C^*$-algebra $C^*(\A)  
\times_{\delta^m} G \times_{\widehat{\delta^m}} H$, hence must be  
contractive, and we arrive at the desired inequality $\< b \cdot x, 
b  
\cdot x \>\rip{C_0} \leq \|b\|^2 \< x,x \>\rip{C_0}$ for $b \in 
B_0$.   
The inequality $\lip{B_0}\< x \cdot c, x \cdot c \> \leq \|c\|^2  
\lip{B_0}\< x,x \>$ is proved similarly. 
\end{proof}

\begin{rem} 
\thmref{imp} can be made equivariant for appropriate  
coactions of $G$ on $C^*(\A)\times_{\delta^m} G 
\times_{\what{\delta^m}} H$ and  $C^*(\c A \times G/H)$.  The 
``appropriate'' coaction of $G$ on $C^*(\A)\times_{\delta^m} G 
\times_{\what{\delta^m}} H$ is  
$\infl(\what{\delta^m}|H)\what{\ }$, the inflation from $H$ to $G$ 
of  
the dual coaction of the restricted action $\what{\delta^m}|H$.  
Recall  
that if $\eps \: B \to B \otimes C^*(H)$ is a coaction of $H$, then  
$\infl \eps = (\id \otimes C_H) \circ \eps \: B \to B \otimes 
C^*(G)$,  
where $C_H \: C^*(H) \to C^*(G)$ denotes the natural inclusion.  The  
``appropriate'' coaction of $G$ on $C^*(\c A \times G/H)$ is the  
\emph{decomposition} 
coaction $\delta\dec$, given on the generators by 
\[ 
\delta\dec(a_s, tH) = (a_s, tH) \otimes s. 
\] 
This is readily verified to give a homomorphism of the Fell bundle 
$\c  
A \times G/H$, hence of the $C^*$-algebra $C^*(\c A \times G/H)$, 
into  
the $C^*$-algebra $C^*(\c A \times G/H) \otimes C^*(G)$.  Moreover 
it  
is easy to check that this homomorphism is nondegenerate and 
satisfies  
the coaction identity, and is injective because $(\id \otimes 1_G)  
\circ \delta\dec$ is the identity map on $C^*(\c A \times G/H)$, 
where  
$1_G$ denotes the trivial one-dimensional representation of $G$. 
 
Recall that $X$ is the completion of the $B_0 - C_0$ 
pre-imprimitivity  
bimodule $X_0$ with operations given by \eqref{eq-X0}.  Adapting 
from  
\cite{er:stabilization} to our context (with non-normal subgroups,  
full coactions, full crossed products, and discrete groups), we find 
a  
unique homomorphism $\delta_X$ of $X_0$ into the $(B_0 \odot C^*(G)) 
-  
(C_0 \odot C^*(G))$ pre-imprimitivity bimodule $X_0 \odot C^*(G)$ 
such  
that 
\[ 
\delta_X(a, t) = (a, t) \otimes t^{-1}. 
\] 
Hence $\delta_X$ extends uniquely to a homomorphism of the 
$C^*(\A)\times_{\delta^m} G \times_{\what{\delta^m}} H - C^*(\c A 
\times G/H)$  
imprimitivity bimodule $X$ into the $\bigl( (C^*(\A)\times_{\delta^m} 
G \times_{\what{\delta^m}} H) \otimes G^*(G) \bigr) - \bigl( C^*(\c 
A  
\times G/H) \otimes G^*(G) \bigr)$ imprimitivity bimodule $X \otimes  
C^*(G)$.  Moreover, it is easy to verify on the generators that this  
homomorphism is nondegenerate and satisfies the coaction identity.   
Therefore $\delta_X$ implements a Morita equivalence between the  
coactions $\infl(\what{\delta^m}|H)\what{\ }$ and $\delta\dec$.  
Observe
that  
in the extreme case $H = G$ we obtain a Morita equivalence between 
$\what{\what{\delta^m}}$ and $\delta^m$, where 
$\what{\what{\delta^m}}$ is the dual coaction 
on the {\em maximal} crossed product 
$C^*(\A)\times_{\delta^m} G \times_{\what{\delta^m}} G$.
Note that this result can not be deduced from the usual Katayama 
duality
theorem, since this only works for normal (or reduced) coactions.
\end{rem} 

\begin{rem}\label{Marke}
Now suppose the subgroup $H$ is normal in $G$.  
Let $\delta^m|$ denote the restriction of the coaction
$\delta^m:C^*(\A)\to C^*(\A)\otimes C^*(G)$ to $G/H$ (see \S2).
 The theory of \cite{eq:induced} shows how to induce $\delta^m|$ to 
a  
coaction $\ind \delta^m|$ of $G$ on an induced $C^*$-algebra $\ind 
C^*(\A)$ 
(we refer to \S 4 below for the precise definitions),  
and \cite[Theorem 4.1]{eq:induced} gives an $\ind C^*(\A) 
\times_{\ind  
\delta^m|} G - C^*(\A) \times_{\delta^m|} G/H$  
imprimitivity bimodule $Z$.  Dually to the situation regarding 
Green's  
imprimitivity theorem for actions, $Z$ is isomorphic to the $C^*(\A)  
\times_{\delta^m} G \times_{\what{\delta^m}} H - 
C^*(\A)\times_{\delta^m|} G/H$  imprimitivity bimodule of 
\thmref{imp}.  The interested reader can  
check that the map 
\[ 
(a_s,t) \mapsto (a_s,st) 
\] 
of $\c A \times G$ extends uniquely to an isomorphism of $X$ onto 
$Z$  
whose left coefficient isomorphism from 
$C^*(\A) \times_{\delta^m} G \times_{\what{\delta^m}} H$ 
onto 
$\ind C^*(\A) \times_{\ind \delta^m|} G$ 
is given on the generators by 
\[ 
(a_s, t, h) \mapsto (a_s, sth^{-1}t^{-1}, th). 
\] 
Note that in order to see that this map on the generators extends to
the desired isomorphism, one {\em uses} the fact that
$X$ {\em is} an $C^*(\A)  
\times_{\delta^m} G \times_{\what{\delta^m}} H - 
C^*(\A)\times_{\delta^m|} G/H$ imprimitivity bimodule, since then
the result follows from an easy application of Lemma \ref{module hom}.
We actually were not able to get this isomorphism directly, i.e., to
deduce Theorem \ref{imp} directly from \cite[Theorem 4.1]{eq:induced}
in the special case where $H$ is normal in $G$.
\end{rem}

In what follows we want to show that the imprimitivity bimodule of 
Theorem  
\ref{imp} factors to an imprimitivity bimodule for the reduced  
cross-sectional algebras.  Thus we want to show 
 
\begin{thm} 
\label{thm-reduced} 
Let $X$ be the $C^*(\A) \times_{\delta^m} G  
\times_{\widehat{\delta^m}} H - C^*(\A \times G/H)$ imprimitivity  
bimodule of Theorem \ref{imp}, and let $I$ and $J$ denote the 
kernels  
of the regular representations of $C^*(\A) \times_{\delta^m} G  
\times_{\widehat{\delta^m}} H$ and $C^*(\A \times G/H)$, 
respectively.   
Then $I$ is equal to the ideal of $C^*(\A) \times_{\delta^m} G  
\times_{\widehat{\delta^m}} H$ induced from $J$ via $X$ and 
therefore  
$Y : = X/(X \cdot J)$ has a canonical structure as a $C^*(\A)  
\times_{\delta^m} G \times_{\widehat{\delta^m},r} H - C_r^*(\A 
\times  
G/H)$ imprimitivity bimodule. 
\end{thm} 
 
Recall that a regular representation of $C^*(\A) \times_{\delta^m} G  
\times_{\widehat{\delta^m}} H$ is by definition an induced  
representation $\Ind_{\{e\}}^H(\pi \times \mu)$, where $\pi \times  
\mu$ is any given faithful representation of $C^*(\A)  
\times_{\delta^m} G$.  Moreover, using Proposition 
\ref{prop-reduced},  
we see that the kernel of the regular representation of $C^*(\A 
\times  
G/H)$ equals $\ker(\pi \times \mu|)$, where $\pi \times \mu|$ 
denotes  
the restriction of $\pi \times \mu$ to $C^*(\A \times G/H)$ (note 
that  
if we represent $C^*(\A) \times_{\delta^m} G$ faithfully on a 
Hilbert  
space $\H$ via $\pi \times \mu$, then we may identify the pair  
$(\pi,\mu)$ with the canonical inclusions $(j_{C^*(\c A)},j_G)$, so  
that Proposition \ref{prop-reduced} applies).  Using these  
realizations of the (kernels of the) regular representations, and 
the  
fact that inducing ideals is compatible with the process of inducing  
representations, that is, $\ker \big(\Ind^X(\pi \times \mu|)\big) =  
\Ind^X \big(\ker(\pi \times \mu|)\big)$, where $X$ is the bimodule 
of Theorem \ref{imp},
Theorem \ref{thm-reduced}  
follows from 
 
\begin{prop} 
\label{prop-indred} 
Let $\pi \times \mu$ be any representation of $C^*(\A)  
\times_{\delta^m} G$.  Then $\Ind_{\{e\}}^H(\pi \times \mu)$ is  
equivalent to $\Ind^X(\pi \times \mu|)$, where $\Ind^X \: 
\Rep(C^*(\A  
\times G/H)) \to \Rep(C^*(\A) \times_{\delta^m} G  
\times_{\widehat{\delta^m}} H)$ denotes induction via the bimodule
$X$ of Theorem \ref{imp}. 
\end{prop} 
 
The above proposition constitutes a (new) version of the duality  
results for induction and restriction of representations obtained in  
\cite{ech:indres, ekr, kqr:resind} 
which works for arbitrary subgroups of discrete groups. 
 
\begin{proof}[Proof of \propref{prop-indred}] 
If $(D, H, \alpha)$ is a system with $H$ discrete, then the process 
of  
inducing representations from $D$ to $D \times_{\alpha} H$ can be  
described as follows: let $W_0 = \Gamma_c(D \times H)$, where $D  
\times H$ denotes the semidirect product bundle.  Then translating  
Green's formulas (see \cite{gre:local}) to this special situation,  
$W_0$ becomes a right pre-Hilbert $D$-module with $D$-valued inner  
product given by 
\[ 
\< (b,h), (c,l) \>\rip{D} = \alpha_{h^{-1}}(b^*c) \nonzero{h = l}. 
\] 
The left convolution action of $\Gamma_c(D \times H)$ then extends 
to  
a $^*$-homomorphism of $D \times_{\alpha} H$ into $\mathcal L_D(W)$,  
where $W$ denotes the Hilbert-module completion of $W_0$.  The map  
$\Ind_{\{e\}}^H \: \Rep(D) \to \Rep(D \times_{\alpha} H)$ is just  
induction via the right-Hilbert $D \times_{\alpha} H - D$ bimodule  
$W$. 
 
In our situation we have $D = C^*(\A) \times_{\delta^m} G$ and 
$\alpha  
= \widehat{\delta^m}|$.  If we restrict the actions and inner 
products  
given above to the dense subalgebras $D_0: = \Gamma_c(\A \times G)$ 
of  
$D$ and $B_0: = \Gamma_c(\A \times G \times H)$ of $B : = D  
\times_{\alpha} H$, and the dense $B_0 - D_0$ subbimodule $Y_0 =  
\Gamma_c(\A \times G \times H)$ of $W$, then $Y_0$ becomes a  
pre-right-Hilbert $B_0 - D_0$ bimodule with inner product and 
actions  
given by the formulas 
\begin{align*} 
\< (a_u,v,l), (b_r,x,m) \>\rip{D_0} 
&= (a_u^* b_r,xl) \nonzero{uv = rx, l = m} 
\\ 
(a_s,t,h) \cdot(b_u,v,l) 
&= (a_s b_u,vh^{-1},hl) \nonzero{th = uv} 
\\ 
(a_u,v,l) \cdot(b_y,z) 
&= (a_u b_y, zl^{-1}, l) \nonzero{vl = yz}. 
\end{align*} 
Clearly, the completion $Y$ of $Y_0$ coincides with $W$ as a  
right-Hilbert $B - D$ bimodule. 
 
We now describe restriction of representations from $C^*(\A)  
\times_{\delta^m} G$ to $C^*(\A \times G/H)$ in bimodule language.   
For this we view $Z : = C^*(\A) \times_{\delta^m} G$ as a  
right-Hilbert $C^*(\A \times G/H) - C^*(\A) \times_{\delta^m} G$  
bimodule equipped with the canonical right $C^*(\A) 
\times_{\delta^m}  
G$-valued inner product and left action of $C^*(\A \times G/H)$ 
given  
by the natural homomorphism of $C^*(\A \times G/H)$ into $M(C^*(\A)  
\times_{\delta^m} G)$.  Then the restriction $\pi \times \mu|$ is  
equivalent to the induced representation $\Ind^Z(\pi \times \mu)$.  
If  
we restrict our attention to the dense subspace $Z_0: = \Gamma_c(\A  
\times G)$ of $Z$, and the dense subalgebras $D_0 = \Gamma_c(\A 
\times  
G)$ and $C_0 = \Gamma_c(\A \times G/H)$ of $D = C^*(\A)  
\times_{\delta^m} G$ and $C = C^*(\A \times G/H)$, respectively, 
then  
the formulas for the inner product and the left action are given on  
the generators by 
\begin{align*} 
\< (a_u,v), (b_r,x) \>\rip{D_0} 
&= (a_u^* b_r, x) \nonzero{uv = rx} 
\\ 
(a_s,tH) \cdot(b_u,v) 
&= (a_s b_u,v) \nonzero{tH = uvH}. 
\end{align*} 
Finally, let $X_0 = \Gamma_c(\A \times G)$ be the $B_0 - C_0$  
pre-imprimitivity bimodule of Theorem \ref{imp} with formulas given 
as  
in Equation \eqref{eq-X0}.  
It is now straightforward to check that the map
$V:X_0\odot Z_0\to Y_0$ defined on the generators by
\begin{equation} 
\label{eq-V} 
V \bigl( (a_s,t) \otimes (a_u,v) \bigr) := (a_s a_u, u^{-1}t, 
t^{-1}uv)  
\nonzero{tH = uvH} 
\end{equation} 
extends to a right-Hilbert $B-D$ bimodule isomorphism of $X\otimes_CZ$
onto $Y$.
This gives the following chain of unitary equivalences among  
representations of $C^*(\c A) \times_{\delta^m} G  
\times_{\what{\delta^m}} H$: 
\begin{align*} 
\Ind_{\{e\}}^H(\pi \times \mu) 
& \cong \Ind^Y(\pi \times \mu) 
\cong \Ind^{X \times_CZ}(\pi \times \mu) 
\\& \cong \Ind^X \big(\Ind^Z(\pi \times \mu)\big) 
\cong \Ind^X(\pi \times \mu|). 
\end{align*} 
\end{proof} 
 
As an immediate corollary of Theorem \ref{thm-reduced} we get the  
following application to amenability of Fell bundles in the sense of  
Exel and Kumjian (see the discussion in \secref{prelim}). 
 
\begin{cor} 
\label{cor-amenable} 
Let $\A$ be a Fell bundle over the discrete group $G$ and let $H$ be 
a  
subgroup of $G$.  Then the following are equivalent: 
\begin{enumerate} 
\item 
$\A \times G/H$ is an amenable Fell bundle over the groupoid $G 
\times  
G/H$. 
 
\item 
$\A \times G \times H$ is an amenable Fell bundle over the groupoid 
$G  
\times G \times H$. 
 
\item 
The semi direct product bundle $(C^*(\A) \times_{\delta^m} G) \times  
H$ is an amenable Fell bundle over $H$. 
\end{enumerate}  
In particular, $\A$ is amenable if and only if the double dual 
bundle  
$(C^*(\A) \times_{\delta^m} G) \times G$ \textup(or $\A \times G  
\times G$\textup) is amenable. 
\end{cor} 
 
Note that in the formulation of the above corollary we could have  
replaced $(C^*(\A) \times_{\delta^m} G) \times H$ with $(A  
\times_{\delta} G) \times H$, where $\delta \: A \to A \otimes 
C^*(G)$  
is any coaction which has $\A$ as an underlying Fell 
bundle (see  
\cite[Lemma 2.1]{eq:induced}).
Just to prevent any misunderstandings: the
 groupoid structure on
$G\times G\times H$ is given by
$$(s,t,h)(u,v,l)=(su,vh^{-1},hl)\quad\text{if}\quad th=uv,$$
which is compatible with Equation (\ref{eq-b0}).
 
We are now going to derive from \thmref{imp} an imprimitivity theorem 
which works for {\em any} coaction of a discrete group $G$.  In order 
to 
prepare the statement recall from  \cite{eq:induced}
 that if $\delta \: A \to A \otimes 
C^*(G)$ is a coaction of $G$ on $A$, and if $(\A, G)$ denotes the 
corresponding Fell bundle, then $A$ is a completion of $\Gamma_c(\A)$ 
with respect to a $C^*$-norm $\|\cdot\|_{\nu}$ which lies between the 
norms $\|\cdot\|_{\max}$ and $\|\cdot\|_{\min}$ arising from viewing 
$\Gamma_c(\A)$ as a dense subalgebra of $C^*(\A)$ and $C_r^*(\A)$, 
respectively.  If $N$ is a normal subgroup of $G$, we may restrict 
$\delta$ to a coaction $\delta| = (\id_A \otimes q) \circ \delta \: 
A \to A \otimes C^*(G/N)$, where $q \: C^*(G) \to C^*(G/N)$ denotes 
the 
quotient map.  If $\delta$ is normal (that is, $A = C_r^*(\A)$), then 
it 
follows from the generalization of Mansfield's imprimitivity theorem 
to 
nonamenable groups obtained in \cite[Corollary 
3.4]{kq:imprimitivity}, 
that $A \times_{\delta|} G/N$ is Morita equivalent to the reduced 
double 
crossed product $A \times_{\delta} G \times_{\widehat{\delta},r} N$, 
but it was not clear at all whether there is a similar imprimitivity 
theorem for general coactions.  For discrete $G$, the following 
theorem 
gives a complete answer to this open question. 
 
\begin{thm} 
\label{general Mans} 
Let $(A, G, \delta)$ be a coaction of the discrete group $G$, and 
let  
$N$ be a normal subgroup of $G$.  Let $(\A,G)$ denote the  
corresponding Fell bundle, and let $\|\cdot\|_{\nu}$ denote the  
$C^*$-norm on $\Gamma_c(\A)$ corresponding to $A$.  Then there 
exists  
a $C^*$-norm $\|\cdot\|_{\mu}$ \textup(lying between  
$\|\cdot\|_{\max}$ and $\|\cdot\|_{\min}$\textup) on $\Gamma_c((A  
\times_{\delta} G) \times N)$ and a quotient $Z$ of the bimodule $X$  
of Theorem \ref{imp} such that $Z$ becomes an 
\[ 
A \times_{\delta} G \times_{\widehat{\delta},\mu} N - 
A \times_{\delta|} G/N 
\] 
imprimitivity bimodule, where $A \times_{\delta} G  
\times_{\widehat{\delta},\mu} N$ denotes the completion of  
$\Gamma_c((A \times_{\delta} G) \times N)$ with respect to  
$\|\cdot\|_{\mu}$.  Moreover, if $\|\cdot\|_{\nu}$ is the minimal  
\textup(respectively, maximal\textup) norm, then $\|\cdot\|_{\mu}$ 
is  
also the minimal \textup(respectively, maximal\textup) norm. 
\end{thm} 
 
\begin{proof} 
Let $\delta^m|$ and $\delta^n|$ denote the restrictions of the  
coactions $\delta^m$ and $\delta^n$ of $G$ on $C^*(\A)$ and  
$C_r^*(\A)$ to $G/N$, respectively, and let 
\[ 
\xymatrix{ 
{C^*(\A)} 
\ar[r]^{\phi} 
\ar[dr]_{\Lambda} 
& {A} 
\ar[d]^{\lambda} 
\\ 
& {C_r^*(\A)} 
} 
\] 
be the commutative diagram of surjections determined by the identity  
map on $\Gamma_c(\A)$.  Since $\phi$, $\lambda$, and $\Lambda$ are  
equivariant with respect to the coactions $\delta^m, \delta$, and  
$\delta^n$, and hence also with respect to their restrictions  
$\delta^m|$, $\delta|$, and $\delta^n|$, we obtain a commutative  
diagram of surjections 
\[ 
\xymatrix{ 
{C^*(\A) \times_{\delta^m|} G/N} 
\ar[r]^-{\phi \times G/N} 
\ar[dr]_{\Lambda \times G/N} 
& {A \times_{\delta|} G/N} 
\ar[d]^{\lambda \times G/N} 
\\ 
& {C_r^*(\A) \times_{\delta^n|} G/N} 
}. 
\] 
Note that $\phi \times G/N$, $\lambda \times G/N$, and $\Lambda 
\times  
G/N$ are all given by the identity on $\Gamma_c(\A \times G/N)$,  
sitting as a dense subalgebra in all three crossed products.   
Moreover, it follows from Corollaries \ref{cor-max} and 
\ref{cor-min}  
that the identity map on $\Gamma_c(\A \times G/N)$ also induces  
isomorphisms $C^*(\A) \times_{\delta^m|} G/N \cong C^*(\A \times 
G/N)$  
and $C_r^*(\A) \times_{\delta^n|} G/N \cong C_r^*(\A \times G/N)$. 
 
Let $X$ be the $C^*(\A) \times_{\delta^m} G  
\times_{\widehat{\delta^m}} N - C^*(\A) \times_{\delta^m|} G/N$  
imprimitivity bimodule of Theorem \ref{imp}.  Since the crossed  
product $A \times_{\delta} G$ and the dual action $\widehat{\delta}$  
only depend on the Fell bundle $\A$ and \emph{not} on the particular  
choice of the cross-sectional algebra $A$ (see \cite[Lemma  
2.1]{eq:induced}) we can replace $C^*(\A) \times_{\delta^m} G  
\times_{\widehat{\delta^m}} N$ on the left hand side with either $A  
\times_{\delta} G \times_{\widehat{\delta}} N$ or $C_r^*(\A)  
\times_{\delta^n} G \times_{\widehat{\delta^n}} N$.  If  
$\|\cdot\|_{\nu} = \|\cdot\|_{\max}$, we have $A = C^*(\A)$ and  
$\delta = \delta^m$, so $X$ is actually an $A \times_{\delta} G  
\times_{\widehat\delta} N - A \times_{\delta|} G/N$ imprimitivity  
bimodule.  Thus we get $\|\cdot\|_{\mu} = \|\cdot\|_{\max}$ in this  
case.  If $\|\cdot\|_{\nu} = \|\cdot\|_{\min}$, (that is, $A =  
C_r^*(\A)$ and $\delta = \delta^n$), then Theorem \ref{thm-reduced}  
shows  that
$X$ factors through an $A \times_{\delta} G  
\times_{\widehat\delta,r} N - A \times_{\delta|} G/N$ imprimitivity  
bimodule, and we get $\|\cdot\|_{\mu} = \|\cdot\|_{\min}$ on the 
left  
hand side. 
 
In general, it follows from the above considerations of the maps 
$\phi  
\times G/N$, $\lambda \times G/N$ and $\Lambda \times G/N$ that the  
kernel $L$ of $\phi \times G/N \: C^*(\A) \times_{\delta^m|} G/N \to 
A  
\times_{\delta|} G/N$ contains the kernel of the regular  
representation $\Lambda \times G/N$, and hence the ideal $K$ of $A  
\times_{\delta} G \times_{\widehat\delta} N$ induced from $L$ via 
$X$  
contains the kernel of the regular representation of $A  
\times_{\delta} G \times_{\widehat\delta} N$ (which by Theorem  
\ref{thm-reduced} is the ideal induced from $\ker \Lambda \times 
G/N$  
via $X$).  It follows that $(A \times_{\delta} G  
\times_{\widehat\delta} N)/K$ is a completion of $\Gamma_c\bigl( (A  
\times_{\delta} G) \times N \bigr)$ with respect to a norm  
$\|\cdot\|_{\mu}$ which lies between $\|\cdot\|_{\max}$ and  
$\|\cdot\|_{\min}$.  Therefore $Z : = X/(X \cdot L)$ is an $A  
\times_{\delta} G \times_{\widehat{\delta},\mu} N - A 
\times_{\delta|}  
G/N$ imprimitivity bimodule. 
\end{proof}


\section{Applications to induced coactions} 
\label{appl} 
 
Recall from \cite{eq:induced} that if $H$ is a normal subgroup
of the discrete group $G$ and $\eps:D\to D\otimes C^*(G/H)$
is a coaction of the quotient group $G/H$ with underlying
Fell bundle $\D$, then the {\em induced coaction} $\delta:=\Ind\eps$
is defined as the dual coaction on the maximal cross-sectional algebra
$C^*(q^*\D)$, where $q^*\D$ denotes the pull back of $\D$ via the
quotient map $q:G\to G/H$. 
In \cite[Theorem 5.6]{eq:induced} we prove an analogue of a theorem 
of  
Olesen and Pedersen, in which we show that
an action of $G$ is twisted in  
the sense of Green over a normal subgroup $H$ if and only if the 
dual  
coaction is induced 
from a coaction  
of $G/H$.  Of course it is a natural question whether the dual of 
this  
result is also true, namely whether a coaction $\delta$ of $G$ is  
induced from $G/H$ if and only if the dual action $\widehat{\delta}$  
is twisted over $H$.  We are now prepared to show that this is indeed the  
case if we assume the coaction is maximal. 
 
\begin{thm} 
\label{twist} 
Let $\delta$ be a coaction of the discrete group $G$ on the  
$C^*$-algebra $A$, and assume that $A = C^*(\c A)$, where $\c A$ is  
the associated Fell bundle.  Then $\delta$ is induced from a 
quotient  
$G/H$ of $G$ by a normal subgroup $H$ if and only if there exists a  
Green-twist $\tau \: H \to UM(A \times_{\delta} G)$ for the dual  
action $\widehat{\delta}$ of $G$ on $A \times_{\delta} G$. 
\end{thm} 
 
\begin{proof} 
Suppose first that $\delta$ is induced from a coaction $\eps \: D \to 
D \otimes C^*(G/H)$ of $G/H$.  This means that $\c A$ is the 
pull-back 
bundle $\pb \D = \{ (D_{sH}, s) : s \in G \}$, where $\D$ is the Fell 
bundle over $G/H$ corresponding to $\eps$. It follows then from 
\cite[Theorem 4.1]{eq:induced} that $\widehat\delta$ 
is Morita equivalent to the inflated action $\Inf \widehat\eps$. 
Since the trivial homomorphism taking $H$ to the identity of $\UM(D 
\times_{\eps} G/H)$ is a twist for $\Inf \widehat\eps$, and since 
Morita 
equivalence of actions preserves the property of being twisted over 
$H$, 
by \cite[Proposition 2]{ech:twisted}, it follows that 
$\widehat{\delta}$ 
is twisted over $H$. 
 
Let us now assume that $\widehat{\delta}$ is twisted over $H$.  Put 
$E  
= A \times_{\delta} G$, and let us look at the semi-direct product  
bundle $E \times G$ corresponding to $\widehat{\delta}$.  The $E  
\times_{\widehat{\delta}} G - A$ imprimitivity bimodule $X$ of  
\thmref{imp} for the case $H = G$ is the completion of the $B_0 - 
C_0$  
pre-imprimitivity bimodule $X_0$ with operations given by the  
corresponding special case of \eqref{eq-X0}.  For $t \in G$ define  
$X_t = \clsp \{(A_s,t^{-1}): s \in G\} \subseteq X$.  Then the  
following assertions are true: 
\begin{enumerate} 
\item 
$(E,s) \cdot X_t = X_{st}$ and 
$X_s \cdot A_t \subseteq X_{st}$; 
 
\item 
$\< X_s, X_t \>\rip A \subseteq A_{s^{-1}t}$ and 
$\< X_e, X_s \>\rip A = A_s$ 
\end{enumerate} 
for all $s,t \in G$, where all spaces here are to be interpreted as  
the respective closed linear spans. 
 
Since, by assumption, $\widehat{\delta}$ is twisted over $H$, the  
semidirect product bundle $E \times G$ is induced from the twisted  
semidirect product bundle $E \times_HG$ over $G/H$, by \cite[Theorem  
5.6]{eq:induced}.  Thus, by \cite[Theorem 5.1]{eq:induced} there  
exists a map $v \: H \to \UM(E \times G)$ with $v_h \in \UM(E, h)$ 
satisfying $(b,t) v_h = v_{tht^{-1}} (b,t)$ for all $(b,t) \in E 
\times  
G$ and $h \in H$ (see \cite[\S5]{eq:induced} for the 
notation).  
If we manage to construct a map $u \: H \to UM(A)$  
with similar properties, the result will follow from the other  
direction of \cite[Theorem 5.1]{eq:induced}. 
 
Since $v_h$ is canonically identified with an element of $UM(E  
\times_{\what\delta} G)$, it acts as a unitary adjointable operator 
on  
the Hilbert $A$-module $X$, with adjoint $v_h^* = v_h^{-1} =  
v_{h^{-1}}$.  Because $X_s = (E,e) \cdot X_s$ by (i), each $v_h$ 
gives  
rise to an isometry $x_s \mapsto v_h \cdot x_s$ from $X_s$ to 
$X_{hs}$  
with inverse $x_{hs} \mapsto v_{h^{-1}} \cdot x_s$.  We use this 
fact  
to define left and right multiplication of $u_h$ with elements in  
$A_s$: 
\begin{equation} 
\label{eq-mult} 
u_h \< x_e, x_s \>\rip A 
= \< x_e, v_h \cdot x_s \>\rip A \in A_{hs} 
\midtext{and} 
\< x_e, x_s \>\rip A u_h 
= \< x_e, v_{shs^{-1}} \cdot x_s \>\rip A \in A_{sh}. 
\end{equation} 
To see that these formulas determine well-defined isometric maps 
$A_s  
\to A_{hs}$ and $A_s \to A_{sh}$, respectively, we first observe 
that  
for all $x_t \in X_t, x_e \in X_e$ and $x_s \in X_s$, 
\begin{equation} 
\label{eq-funny} 
\begin{split} 
x_t \cdot \< x_e, v_h \cdot x_s \>\rip A 
&= \lip{B_0}\< x_t, x_e \> v_h \cdot x_s 
= v_{tht^{-1}} \lip{B_0}\< x_t, x_e \> \cdot x_s 
\\&= \lip{B_0}\< v_{tht^{-1}} \cdot x_t, x_e \> \cdot x_s 
= v_{tht^{-1}} \cdot x_t \cdot \< x_e, x_s \>\rip A 
\end{split} 
\end{equation} 
Using \eqref{eq-funny} for $t = e$ we compute 
\begin{align*} 
& \biggl\| \sum_{i=1}^n \< x_e^i, v_h \cdot x_s^i \>\rip A \biggr\|^2 
\\&\quad = \biggl\| \sum_{i,j=1}^n \< v_h \cdot x_s^j, x_e^j \>\rip 
A  
\< x_e^i, v_h \cdot x_s^i \>\rip A \biggr\| 
= \biggl\| \sum_{i,j=1}^n \bigl\< v_h \cdot x_s^j, x_e^j \cdot 
\< x_e^i, v_h \cdot x_s^i \>\rip A \bigr\>\rip A \biggr\| 
\\&\quad\overset{\eqref{eq-funny}}{=} 
\biggl\| \sum_{i,j=1}^n \bigl\< v_h \cdot x_s^j, 
v_h \cdot x_e^j \cdot \< x_e^i, x_s^i \>\rip A \bigr\>\rip A \biggr\| 
= \biggl\| \sum_{i,j=1}^n \< v_h \cdot x_s^j, v_h \cdot x_e^j \>\rip 
A  
\< x_e^i, x_s^i \>\rip A \biggr\| 
\\&\quad = \biggl\| \sum_{i,j=1}^n \< x_s^j, x_e^j \>\rip A 
\< x_e^i, x_s^i \>\rip A \biggr\| 
= \biggl\| \sum_{i=1}^n \< x_e^i, x_s^i \>\rip A \biggr\|^2. 
\end{align*} 
Thus, left multiplication with $u_h$, as defined in Equation  
\eqref{eq-mult}, is a well-defined isometry with inverse 
$u_{h^{-1}}$,  
and replacing $h$ by $shs^{-1}$ gives the similar result for right  
multiplication.  We show $(a_su_h)a_t = a_s(u_ha_t)$ for all $s,t 
\in  
G$ and $h \in H$: 
\begin{align*} 
&\bigl( \< x_e, x_s \>\rip A u_h \bigr) \< y_e, y_t \>\rip A 
= \< x_e, v_{shs^{-1}} \cdot x_s \>\rip A \< y_e, y_t \>\rip A 
\\&\quad = \bigl\< x_e, v_{shs^{-1}} \cdot x_s \cdot 
\< y_e, y_t \>\rip A \bigr\>\rip A 
= \bigl\< x_e, \lip{B_0}\< v_{shs^{-1}} \cdot x_s, y_e \> 
\cdot y_t \bigr\>\rip A 
\\&\quad = \bigl\< x_e, v_{shs^{-1}} 
\lip{B_0}\< x_s, y_e \> \cdot y_t \bigr\>\rip A 
= \bigl\< x_e, \lip{B_0}\< x_s, y_e \> v_h \cdot y_t \bigr\>\rip A 
\\&\quad = \< x_e, x_s \>\rip A \< y_e, v_h \cdot y_t \>\rip A 
= \< x_e, x_s \>\rip A \bigl(u_h \< y_e, y_t \>\rip A\bigr). 
\end{align*} 
To see that $u_h \in \UM(A_h)$ it now suffices to check that $u_h^* 
=  
u_{h^{-1}}$ ($ = u_h^{-1}$).  Using \eqref{eq-funny} we first compute 
\begin{align*} 
&\bigl( \< x_r, x_t \>\rip A u_h \bigr) \< x_e, x_s \>\rip A 
= \bigl\< x_r, x_t \cdot \< x_e, v_h \cdot x_s \>\rip A \bigr\>\rip A 
\\&\quad = \bigl\< x_r, v_{tht^{-1}} \cdot x_t \cdot 
\< x_e, x_s \>\rip A \bigr\>\rip A 
= \< x_r, v_{tht^{-1}} \cdot x_t \>\rip A \< x_e, x_s \>\rip A, 
\end{align*} 
from which it follows that right multiplication of $u_h$ with an 
inner  
product $\< x_r,x_t \>\rip A \in A_{r^{-1}t}$ is given by the 
formula  
$\< x_r, x_t \>\rip A u_h = \< x_r, v_{tht^{-1}} \cdot x_t \>\rip 
A$.   
Using this we now compute 
\begin{align*} 
(u_h \< x_e, x_s \>\rip A)^* 
&= \< x_e, v_h \cdot x_s \>\rip A^* 
= \< v_h \cdot x_s, x_e \>\rip A 
\\&= \< x_s, v_{h^{-1}} \cdot x_e \>\rip A 
= \< x_s, x_e \>\rip A u_{h^{-1}} 
= \< x_e,x_s \>\rip A^* u_{h^{-1}}, 
\end{align*} 
which proves $u_h^* = u_{h^{-1}}$.  Since it follows directly from 
the  
definition of left and right multiplication with $u_h$ that  
$u_{shs^{-1}} a_s = a_s u_h$ for all $a_s \in A_s$, we see that $h  
\mapsto u_h \in \UM(A_h)$ satisfies all requirements of 
\cite[Theorem  
5.1]{eq:induced}, and the result follows. 
\end{proof}

 

\providecommand{\bysame}{\leavevmode\hbox to3em{\hrulefill}\thinspace}

\end{document}